\documentclass[12pt]{article}
\usepackage{epic,eepic,epsf,epsfig}
\usepackage{amsmath,enumerate}
\usepackage{amssymb}
\usepackage{amsbsy}

\usepackage{amsfonts}
\newtheorem{theorem}{Theorem}[section]%
\newtheorem{lemma}[theorem]{Lemma}%
\newtheorem{prop}[theorem]{Proposition}%

\newtheorem{remark}[theorem]{Remark}%

 \def\Omega{\Omega}

\def\f{\noindent}

\newcommand{\qed}{\mbox{\raisebox{0.7ex}{\fbox{}}} \vspace{4truemm}}
\def\demo{\f {\bf Proof.}\hskip10pt}

\begin{document}

\baselineskip 16pt

\title{ \vspace{-1.2cm}
On generalizations of Iwasawa's theorem
\thanks{\scriptsize This research was supported in part by Shandong Provincial Natural Science Foundation, China (ZR2017MA022)
and NSFC (11761079).
\newline
 \hspace*{0.5cm} \scriptsize $^{\ast\ast}$Corresponding
  author.
\newline
       \hspace*{0.5cm} \scriptsize{E-mail addresses:}
       shijt2005@pku.org.cn\,(J. Shi),\,\,xufj2023@s.ytu.edu.cn\,(F. Xu),\,\,shanmj2023@s.ytu.edu.cn\,(M. Shan).}}

\author{Jiangtao Shi\,$^{\ast\ast}$,\,\,Fanjie Xu,\,\, Mengjiao Shan\\
\\
{\small School of Mathematics and Information Sciences, Yantai University, Yantai 264005, China}}

\date{ }

\maketitle \vspace{-.8cm}

\begin{abstract}
Iwasawa's theorem indicates that a finite group $G$ is supersolvable if and only if all maximal chains of the identity in $G$ have
the same length. As generalizations of Iwasawa's theorem, we provide some characterizations of the structure of a finite group $G$ in which all maximal
chains of every minimal subgroup have the same length. Moreover, let $\delta(G)$ be the number of subgroups of $G$ all of whose maximal chains in $G$
do not have the same length, we prove that $G$ is a non-solvable group with $\delta(G)\leq 16$ if and only if $G\cong A_5$.

\medskip \f {\bf Keywords:} supersolvable group; maximal chains; Iwasawa's theorem; solvable group; Fitting subgroup
\end{abstract}

\section{Introduction}\label{s1}

In this paper all groups are considered to be finite. Let $G$ be a group and $H$ a subgroup of $G$. A maximal chain of $H$ in $G$ has length
$r$ means a chain of subgroups $H=H_0<H_1<H_2<\ldots<H_r=G$ such that $H_i$ is maximal in $H_{i+1}$ for every $1\leq i\leq r-1$.

Since every maximal subgroup of a supersolvable group $G$ has prime index, all maximal chains of every subgroup in $G$ have the same length.
Moreover, the following proposition obviously holds for maximal chains of subgroups.

\begin{prop}\ \ \label{p1} Let $G$ be a group, $H<K<G$. If all maximal chains of $H$ in $G$ have the same length, then

$(1)$ all maximal chains of $H$ in $K$ have the same length;

$(2)$ all maximal chains of $K$ in $G$ have the same length.
\end{prop}

\begin{remark}\ \ \label{r1} {\rm Let $G$ be a group, $H<K<G$. If all maximal chains of $K$ in $G$ have the same length, we cannnot ensure that all
maximal chains of $H$ in $G$ have the same length. For example, let $G =A_4=\langle(123),\,(124)\rangle$, $K=\langle(123)\rangle$ and $H=1$. Since $K$
is maximal in $G$, all maximal chains of $K$ in $G$ have the same length. For $H$, it is easy to see that $H<K<G$ and $H<\langle(12)(34)\rangle<\langle(12)(34),\, (13)(24)\rangle<G$
are two maximal chains of $H$ in $G$ having different length.}
\end{remark}

By Proposition~\ref{p1}, if all maximal chains of the identity in $G$ have the same length, then all maximal chains of every subgroup
in $G$ have the same length. Using maximal chains of subgroups Iwasawa provided the following equivalent characterization of supersolvable groups.

\begin{theorem}\ \ \label{thx} {\rm\cite[Theorem 10.3.5]{rob}} A group $G$ is supersolvable if and only if all maximal chains of the identity in $G$ have the same length.
\end{theorem}

As a generalization of Iwasawa's theorem, considering maximal chains of some particular minimal subgroups, we have the following result whose proof is given in Section~\ref{s3}.

\begin{theorem}\ \ \label{l1} Let $G$ be a group. If all maximal chains of every minimal subgroup of order both $2$ and $3$ in $G$ have the same length, then $G$ is solvable.
\end{theorem}

Furthermore, considering maximal chains of every minimal subgroup, we obtain the following result, the proof of which is given in Section~\ref{s4}.

\begin{theorem}\ \ \label{th1} Let $G$ be a group and $F(G)$ the Fitting subgroup of $G$. If all maximal chains of every minimal subgroup in $G$ have the same length, then

$(1)$ $G$ is non-supersolvable if and only if $G=F(G)\rtimes M$ and $\it\Phi$$(G)=Z(G)=1$, where $M$ acts faithfully on $F(G)$, $F(G)$ is the unique minimal normal
subgroup of $G$ and $M$ is a supersolvable maximal subgroup of $G$, $|F(G)|=p^n,\,n\geq 2$;

$(2)$ If $G$ is non-supersolvable, then $F(G)$ is a Sylow subgroup of $G'$. In particular, if $Z(G')\neq 1$, then $F(G)=G'$.
\end{theorem}

\begin{remark}\ \ \label{r2} {\rm If we assume that all maximal chains of every subgroup of $G$ of order a square of a prime in $G$ have the same length,
we cannot ensure that $G$ is solvable. For example, let $G=A_5$. It is easy to see that all maximal chains of subgroup of $A_5$ of order 4 in $A_5$ have the
same length, but $A_5$ is non-solvable.}
\end{remark}

\begin{remark}\ \ \label{r3} {\rm Assume that $G$ is a solvable group, we cannot ensure that all maximal chains of every minimal subgroup in $G$ have
the same length. For example, let $G=S_4=\langle(12),\,(13),\,(14)\rangle$. Take $H=\langle (12)\rangle$ being a minimal subgroup of $G$. It is
easy to see that $H<S_3=\langle(12),\,(13)\rangle<G$ and $H<\langle(12),\,(34)\rangle<D_8=\langle(1423),\,(34)\rangle<G$ are two maximal chains
of $H$ in $G$ having different length.}
\end{remark}

\begin{remark}\ \ \label{r4} {\rm The group in Theorem~\ref{th1} might not be supersolvable. For example, all maximal chains of every
minimal subgroup in $A_4$ have the same length, but $A_4$ is non-supersolvable.}
\end{remark}

\begin{remark}\ \ \label{r5} {\rm In Theorem~\ref{th1} even if assume that $G$ has odd order, we cannot ensure that $G$ is supersolvable.
For example, let $G=\langle a,\,b,\,c\mid\,a^5=b^5=c^3=1,\,[a,b]=1,\,c^{-1}ac=(ab)^{-1},\,c^{-1}bc=a\rangle$. By the definition of $G$,
it is easy to see that all maximal chains of every minimal subgroup in $G$ have the same length but $G$ is non-supersolvable.}
\end{remark}

Note that two examples in Remarks~\ref{r3} and ~\ref{r4} also indicate that even if all maximal chains of every minimal subgroup in $G$ have the
same length, all maximal chains of the identity in $G$ might not have the same length.

\begin{remark}\ \ \label{r6} {\rm Suppose that $G$ is a group satisfying all hypotheses in Theorem~\ref{th1} (2), we cannot ensure that all maximal
chains of every minimal subgroup in $G$ have the same length. For example, let $G=S_4=\langle (12),\,(13),\,(14)\rangle$, one has $F(G)=K_4=\langle(12)(34),\,(13)(24)\rangle$.
Take $M=S_3=\langle (12),\,(13)\rangle$, then $G=F(G)\rtimes M$. Note that $G$ satisfies all hypotheses in Theorem~\ref{th1} (2), but all maximal chains of
$\langle (12)\rangle$ in $G$ have different length.}
\end{remark}

Denote by $\delta(G)$ the number of subgroups of $G$ all of whose maximal chains in $G$ do not have the same length. Observe that $\delta(A_5)=16$ including
the identity and 15 subgroups of order 2. In this paper, we obtain the following result whose proof is given in Section~\ref{s5}.

\begin{theorem}\ \ \label{th2} Let $G$ be a non-solvable group. Then $\delta(G)\leq 16$ if and only if $G\cong A_5$.
\end{theorem}

In Section~\ref{s6}, we also discuss the cases that all maximal chains of some other special subgroups in $G$ have the same length.

\section{Some necessary lemmas}\label{s2}

\begin{lemma}\ \ \label{l01} Let $G$ be a group. Then $G$ is supersolvable if and only if all maximal chains of $\it\Phi$$(G)$ in $G$ have the same length.
\end{lemma}

\demo We only need to prove the sufficiency part. Since all maximal chains of $\it\Phi$$(G)$ in $G$ have the same length, all maximal chains of the
identity subgroup $\bar 1=\it\Phi$$(G)/\it\Phi$$(G)$ in $G/\it\Phi$$(G)$ have the same length. By {\rm\cite[Theorem 10.3.5]{rob}}, $G/\it\Phi$$(G)$
is supersolvable. It follows that $G$ is supersolvable. \hfill\qed

\begin{lemma} {\rm\cite{Thompson}}\ \ \label{l02} A non-abelian simple group
all of whose proper subgroups are solvable is said to be a minimal
simple group, there are five classes in all:

$(1)$ $PSL_2(p)$, where $p>3$ and $5\nmid p^2-1$;

$(2)$ $PSL_2(2^q)$, where $q$ is a prime;

$(3)$ $PSL_2(3^q)$, where $q$ is an odd prime;

$(4)$ $PSL_3(3)$;

$(5)$ $S_z(2^q)$, where $q$ is an odd prime.
\end{lemma}

\begin{lemma} {\rm\cite[Theorem 5.4]{chen}}\ \ \label{l03} Let $G$ be a group of order $n$. If $(n,\,15)=1$, then $G$ is solvable.
\end{lemma}

\begin{lemma} {\rm\cite[Theorem 10.4.2]{rob}}\ \ \label{l04} Let $G$ be a group having a nilpotent maximal subgroup of odd order.
Then $G$ is solvable.
\end{lemma}

\section{Proof of Theorem~\ref{l1}}\label{s3}

\demo Let $G$ be a counterexample of minimal order.

Note that if all maximal chains of every minimal subgroup in $G$ have the same length, then all maximal chains of every non-trivial subgroup
in $G$ have the same length by Proposition~\ref{p1}.

For any maximal subgroup $M$ of $G$, it is obvious that all maximal chains of every minimal
subgroup of order both 2 and 3 in $M$ have the same length by the hypothesis (Note that if there exists a maximal subgroup $E$ of $G$ such that $2\nmid|E|$ or $3\nmid|E|$,
then $E$ naturally satisfies that all maximal chains of every minimal
subgroup of order both 2 and 3 in $E$ have the same length.). By the minimality of $G$, $M$ is solvable. It follows that $G$ is a minimal non-solvable group and then $G/\it\Phi$$(G)$
is a minimal non-abelian simple group.

Claim $\it\Phi$$(G)=1$. Otherwise, assume $\it\Phi$$(G)\neq1$. If $2\mid|\it\Phi$$(G)|$ or $3\mid|\it\Phi$$(G)|$, one has that all maximal chains of $\it\Phi$$(G)$ in $G$ have the same length.
By Lemma~\ref{l01}, $G$ is supersolvable, a contradiction. If $(6,|\it\Phi$$(G)|)=1$, then all maximal chains of every minimal subgroup of $G/\it\Phi$$(G)$ of order both 2 and
3 in $G/\it\Phi$$(G)$ have the same length. It follows that $G/\it\Phi$$(G)$ is solvable by the minimality of $G$, which implies that $G$ is solvable, also a contradiction.

Therefore, $\it\Phi$$(G)=1$. One has that $G$ is a minimal non-abelian simple group. By Lemma~\ref{l02}, we divide our analyses into the following
five cases.

$(1)$ Let $G=PSL_2(p)$, where $p>3$, $(5,\,p^2-1)=1$ and $|G|=\frac{p(p^2-1)}{2}$.

Since $G$ is non-solvable and $(5,\,p^2-1)=1$, one has $3\mid\,p^2-1$ by
Lemma~\ref{l03}.

First, suppose $p^2\equiv\,1\,(\rm mod\,16)$. Then $G$ has a maximal subgroup $M\cong S_4$ by {\rm\cite{DIC}}. It is easy to see that all maximal chains of
subgroups of order 2 in $S_4$ do not have the same length, which implies that all maximal chains of subgroups of order 2 in $G$ do not have the same
length, a contradiction.

Next, suppose $p^2\not\equiv\,1\,(\rm mod\, 16)$. Then $G$ has a maximal subgroup $M\cong A_4$ by {\rm\cite{DIC}}. Consider a subgroup $N$ of $M$ of order 3,
one has that $N<M<G$ is a maximal chain of $N$ in $G$. Let $P$ be a Sylow 3-subgroup of $G$ such that $N\leq P$. By Lemma~\ref{l04}, $P$ is not a
maximal subgroup of $G$. Assume $|P|=3^n$, where $n\geq 1$.

If $n\geq2$, it is easy to see that all maximal chains of $N$ in $G$ do not have the same length, a contradiction. Therefore, $N=P$ is a Sylow 3-subgroup
of $G$.

Since $3\mid\,p^2-1$, one has $3\mid\, p+1$ or $3\mid\,p-1$. It follows that $3\mid\,\frac{p+1}{2}$ or $3\mid\,\frac{p-1}{2}$.

If $3\mid\,\frac{p+1}{2}$,
let $L$ be a maximal subgroup of $G$ that is isomorphic to a dihedral group of order $2\cdot\frac{p+1}{2}$. By Sylow theorem, we can assume $N<L$. Since
all maximal chains of $N$ in $G$ have the same length, one must have $\frac{p+1}{2}=3$, which implies $p=5$. Then $G=PSL_2(5)$. It is easy to see that
all maximal chains of subgroups of order 2 in $PSL_2(5)$ do not have the same length, a contradiction.

If $3\mid\,\frac{p-1}{2}$, let $R$ be a maximal
subgroup of $G$ that is a dihedral group of order $2\cdot\frac{p-1}{2}$. By Sylow theorem, we can assume $N<R$. Since all maximal chains of $N$ in $G$
have the same length, one must have $\frac{p-1}{2}=3$, which implies $p=7$, this contradicts $p^2\not\equiv\,1\,(\rm mod\, 16)$.

$(2)$ Let $G=PSL_2(2^q)$, where $q$ is a prime.

By {\rm\cite{DIC}}, it is easy to see that $G$ has two distinct maximal subgroups $M_1$ and $M_2$ such that
$M_1$ is a dihedral group of order $2\cdot(2^q-1)$ and $M_2$ is the normalizer of the Sylow 2-subgroup of $G$ of order $2^q\cdot(2^q-1)$, where subgroups
of $M_1$ and $M_2$ of order $2^q-1$ are cyclic. Let $N_1$ be a subgroup of $M_1$ of order 2, and let $P$ be a Sylow 2-subgroup of $G$ such that $P<M_2$.
By Sylow theorem, we can assume $N_1\leq P$. Since $q\geq 2$, it is easy to see that all maximal chains of $N_1$ in $G$ do not have the same length,
a contradiction.

$(3)$ Let $G=PSL_2(3^q)$, where $q$ is an odd prime.

By {\rm\cite{DIC}}, $G$ has the following two distinct maximal subgroups: $M_1\cong A_4$ and $M_2$ is the
normalizer of the Sylow 3-subgroup of $G$ of order $3^q\cdot\frac{(3^q-1)}{2}$,  where subgroups of $M_2$ of order $\frac{3^q-1}{2}$ are cyclic. Let $N_1$
be a subgroup of $M_1$ of order 3 and let $Q$ be a Sylow 3-subgroup of $G$ such that $Q<M_2$. By Sylow theorem, we can assume $N_1\leq Q$. Since $q\geq 3$,
it is easy to see that all maximal chains of $N_1$ in $G$ do not have the same length, a contradiction.

$(4)$ Let $G=PSL_3(3)$.

By {\rm\cite{DIC}}, $G$ has a maximal subgroup $M\cong S_4$. It is obvious that there exists a subgroup $S$ of $M$ of order 2
such that all maximal chains of $S$ in $M$ do not have the same length, which implies that all maximal chains of $S$ in $G$ do not have the same length,
a contradiction.

$(5)$ Let $G=S_z(2^q)$, where $q$ is an odd prime.

By {\rm\cite{SUZ}}, $G$ has two distinct maximal subgroups: $M_1$ is the normalizer of the Sylow 2-subgroup
of order $2^{2q}\cdot(2^q-1)$, $M_2$ is a dihedral group of order $2\cdot(2^q-1)$,  where subgroups of $M_1$ and $M_2$ of order $2^q-1$ are cyclic.
Let $N_2$ be a subgroup of $M_2$ of order 2 and let $P$ be a Sylow 2-subgroup of $G$ such that $P<M_1$. By Sylow theorem, we can assume $N_2\leq P$.
Note that $q\geq 3$ and then $2q\geq 6$. It is easy to see that all maximal chains of $N_2$ in $G$ do not have the same length, a contradiction.

By above arguments, the counterexample of minimal order does not exist and so $G$ is solvable.\hfill\qed

\section{Proof of Theorem~\ref{th1}}\label{s4}

\begin{lemma}\ \ \label{l2} Let $G$ be a group in which all maximal chains of every minimal subgroup have the same length, $F(G)$ be a Fitting
subgroup of $G$. Then $G$ is non-supersolvable if and only if $G=F(G)\rtimes M$ and $\it\Phi$$(G)=Z(G)=1$, where $M$ acts faithfully on $F(G)$,
$F(G)$ is the unique minimal normal subgroup of $G$ and $M$ is a supersolvable maximal subgroup of $G$, $|F(G)|=p^n$, $n\geq 2$.
\end{lemma}

\demo Since any minimal normal subgroup of a supersolvable group has prime order, the sufficiency part holds.

In the following we prove the necessity part.

Claim $\it\Phi$$(G)=1$. Otherwise, assume $\it\Phi$$(G)\neq1$. Then all maximal chains of $\it\Phi$$(G)$ in $G$ have the same length. One has that $G$
is supersolvable by Lemma~\ref{l01}, a contradiction. Hence $\it\Phi$$(G)=1$.

Claim $Z(G)=1$. Otherwise, assume $Z(G)\neq1$. Then all maximal chains of $Z(G)$ in $G$ have the same length. For the quotient group $G/Z(G)$, since
all maximal chains of the identity subgroup $Z(G)/Z(G)$ in $G/Z(G)$ have the same length, $G/Z(G)$ is supersolvable by {\rm\cite[Theorem 10.3.5]{rob}},
which implies that $G$ is supersolvable, a contradiction. Therefore, $Z(G)=1$.

Note that $G$ is solvable by Theorem~\ref{l1}. Let $N$ be a minimal normal subgroup of $G$, where $|N|=p^n$, $n\geq1$. Then all maximal chains of $N$ in $G$ have
the same length. Consider the quotient group $G/N$. Since all maximal chains of the identity subgroup $N/N$ in $G/N$ have the same length, $G/N$ is
supersolvable by {\rm\cite[Theorem 10.3.5]{rob}}. If $G$ has at least two distinct minimal normal subgroups, let $N_0$ be another minimal normal subgroup
of $G$ such that $N_0\neq N$, then $N\cap N_0=1$. Arguing as above, one has that $G/N_0$ is supersolvable. It follows that $G\cong
G/(N\cap N_0)$ is isomorphic to a subgroup of $G/N\times G/N_0$ and then $G$ is supersolvable, a contradiction. Therefore, $N$ is the unique minimal normal
subgroup of $G$.

It follows that $F(G)=N$ since $\it\Phi$$(G)=1$ by {\rm\cite[Chapter III, Theorem 4.5]{huppert}}. So $F(G)\nleq\it\Phi$$(G)$. There exists a maximal
subgroup $M$ of $G$ such that $F(G)\nleq M$. One has $G=F(G)M$. Since $G$ is solvable, $F(G)$ is an elementary abelian group. One has $F(G)\cap M\unlhd F(G)$.
Moreover, as $F(G)\cap M\unlhd M$, one has $F(G)\cap M\unlhd G$. Note that $F(G)\cap M<F(G)$, then $F(G)\cap M=1$ by the minimality of $F(G)$. Therefore,
$G=F(G)\rtimes M$. It implies that $M\cong G/F(G)$ is supersolvable.

Since $G$ is solvable, one has $C_G(F(G))\leq F(G)$ by {\rm\cite[Chapter III, Theorem 4.2]{huppert}}. Moreover, as $F(G)\leq C_G(F(G))$, it follows that
$C_G(F(G))=F(G)$, which implies that $M$ acts faithfully on $F(G)$.

Note that $G/F(G)$ is supersolvable and $G$ is non-supersolvable. Therefore, $|F(G)|=p^n$, $n\geq 2$.\hfill\qed

\begin{lemma}\ \ \label{l3} Let $G$ be a group in which all maximal chains of every minimal subgroup have the same length. If $G$ is non-supersolvable,
then $F(G)$ is a Sylow subgroup of $G'$. In particular, if $Z(G')\neq 1$, then $F(G)=G'$.
\end{lemma}

\demo Since $G$ is non-supersolvable, one has that $G/F(G)$ is supersolvable by Lemma~\ref{l2}, where $|F(G)|=p^n$, $n\geq 2$. By {\rm\cite[Theorem 10.5.4]{hall}},
$(G/F(G))'=G'F(G)/F(G)$ is nilpotent. It is obvious that $G'\neq 1$. Since $F(G)$ is the unique minimal normal subgroup of $G$, $F(G)\leq G'$, which implies that
$G'/F(G)$ is nilpotent.

Let $P\in{\rm Syl}_p(G')$. Then $F(G)\leq P$. If $F(G)<P$, then $1<P/F(G)\,{\rm char}\, G'/F(G)\unlhd G/F(G)$ since $G'/F(G)$ is nilpotent.
One has $P\unlhd G$, which implies that $P\leq F(G)$, a contradiction. Hence $F(G)=P\in{\rm Syl}_p(G')$.

Assume $Z(G')\neq 1$. Since $Z(G')\,{\rm char}\, G'\unlhd G$, $Z(G')\unlhd G$. One has $Z(G')=F(G)$. Then $G'/Z(G')=G'/F(G)$ is nilpotent. It follows that
$G'$ is nilpotent and then $G'\leq F(G)$. Moreover, as $F(G)\leq G'$, one has $F(G)=G'$.\hfill\qed

{\bf Proof of Theorem~\ref{th1}.} Combine Lemmas~\ref{l2} and ~\ref{l3} together, we complete the proof of Theorem~\ref{th1}.\hfill\qed

\section{Proof of Theorem~\ref{th2}}\label{s5}

\begin{lemma}\ \ \label{l21} Let $G$ be a group. If $\delta(G)<16$, then $G$ is solvable.
\end{lemma}

\demo Let $G$ be a counterexample of minimal order.

For any proper subgroup $M$ of $G$, if there exists a subgroup $H$ of $M$ such that all maximal chains of $H$
in $M$ do not have the same length, then all maximal chains of $H$ in $G$ do not have the same length. Therefore, $\delta(M)\leq\delta(G)<16$. By the minimality
of $G$, $M$ is solvable. It follows that $G$ is a minimal non-solvable group and then $G/\it\Phi$$(G)$ is a minimal non-abelian simple group.

If $\it\Phi$$(G)\neq 1$, since $\delta(G/\it\Phi$$(G))\leq\delta(G)<16$, $G/\it\Phi$$(G)$ is solvable by the minimality of $G$. It implies that $G$ is solvable, a contradiction.
Therefore, $\it\Phi$$(G)=1$ and then $G$ is a minimal non-abelian simple group.

Note that if all maximal chains of the subgroup $H$ in $G$ do not have the same length, then for any $g\in G$, all maximal chains of $H^g$ in $G$ do not have the same length, too.

By Lemma~\ref{l02}, we divide our analyses into the following five cases.

$(1)$ Let $G=PSL_2(p)$, where $p>3$, $(5,\,p^2-1)=1$ and $|G|=\frac{p(p^2-1)}{2}$.

If $p^2\not\equiv\,1\,(\rm mod\, 16)$. Arguing as in proof of Theorem~\ref{l1}, one has $3\mid\,\frac{p+1}{2}$ or $3\mid\,\frac{p-1}{2}$.

Case$(i)$: Assume $3\mid\,\frac{p+1}{2}$. When $p=5$, then $G=PSL_2(5)\cong A_5$. However, $\delta(A_5)=16$, a contradiction. When $p\geq 7$, $G$ has a subgroup $H$ of order
3 all of whose maximal chains in $G$ do not have the same length. Observe that $|N_G(H)|=2\cdot\frac{p+1}{2}=p+1$. Then $\delta(G)\geq|G:N_G(H)|=\frac{\frac{p(p^2-1)}{2}}{p+1}
=\frac{p(p-1)}{2}\geq 21>16$, a contradiction.

Case $(ii)$: Assume $3\mid\,\frac{p-1}{2}$. As $p^2\not\equiv\,1\,(\rm mod\, 16)$, one has $p>7$. Then $G$ has a subgroup $K$ of order 3 all of whose maximal chains in $G$ do
not have the same length. Observe that $|N_G(K)|=2\cdot\frac{p-1}{2}=p-1$. Then $\delta(G)\geq|G:N_G(K)|=\frac{\frac{p(p^2-1)}{2}}{p-1}=\frac{p(p+1)}{2}>28>16$, also a contradiction.

Next assume $p^2\equiv\,1\,(\rm mod\, 16)$. Then $p\geq 7$.

If $p=7$, then $G=PSL_2(7)$. It is easy to see that $G$ has a subgroup $H$ of order 2 all of whose maximal chains in $G$ do not have the same length. By {\rm\cite{con}}, all
subgroups of $PSL_2(7)$ of order 2 are conjugate in $G$. Then $\delta(G)\geq |G:N_G(H)|=21>16$, a contradiction.

If $p>7$, then $p\geq 17$. By {\rm\cite{DIC}}, $G$ has the following four distinct maximal subgroups: $M_1\cong S_4$, $M_2$ is a dihedral group of order $2\cdot\frac{p+1}{2}$,
$M_3$ is a dihedral group of order $2\cdot\frac{p-1}{2}$ and $M_4$ is the normalizer of the Sylow $p$-subgroup of $G$ of order $\frac{p(p-1)}{2}$. It is easy to see that $G$
has a subgroup $H$ of order 2 all of whose maximal chains in $G$ do not have the same length, and only one of the following three cases might be true by the hypothesis:
$(a)$ $|N_G(H)|\leq|S_4|=24$; $(b)$ $|N_G(H)|\leq2\cdot\frac{p+1}{2}=p+1$; $(c)$ $|N_G(H)|\leq2\cdot\frac{p-1}{2}=p-1$. Note that $p\geq 17$. For case $(a)$, one has
$\delta(G)\geq |G:N_G(H)|\geq\frac{\frac{p(p^2-1)}{2}}{24}\geq 102>16$,
a contradiction. For case $(b)$, one has $\delta(G)\geq |G:N_G(H)|\geq \frac{\frac{p(p^2-1)}{2}}{p+1}=\frac{p(p-1)}{2}\geq 136>16$, a contradiction. For case $(c)$, one has
$\delta(G)\geq |G:N_G(H)|\geq\frac{\frac{p(p^2-1)}{2}}{p-1}=\frac{p(p+1)}{2}\geq 153>16$, also a contradiction.

$(2)$ Let $G=PSL_2(2^q)$, where $q$ is a prime and $|G|=2^q(2^{2q}-1)$.

Arguing as in proof of Theorem~\ref{l1}, $G$ has a subgroup $H$ of order 2 all of whose maximal chains in $G$ do
 not have the same length and $|N_G(H)|=2^q$. Then $\delta(G)\geq |G:N_G(H)|+1=\frac{2^q(2^{2q}-1)}{2^q}+1=2^{2q}-1+1\geq 16$, a contradiction.

$(3)$ Let $G=PSL_2(3^q)$, where $q$ is an odd prime number and $|G|=\frac{3^q(3^{2q}-1)}{2}$.

Arguing as in proof of Theorem~\ref{l1}, $G$ has a subgroup $H$ of order 3 all of whose maximal chains in $G$ do
 not have the same length and $|N_G(H)|\leq\frac{3^q(3^q-1)}{2}$. Then $\delta(G)\geq
|G:N_G(H)|\geq\frac{(\frac{3^q(3^{2q}-1)}{2})}{(\frac{3^q(3^q-1)}{2})}=3^q+1\geq28>16$, a contradiction.

$(4)$ Let $G=PSL_3(3)$.

Arguing as in proof of Theorem~\ref{l1}, $G$ has a subgroup $H$ of order 2 all of whose maximal chains in $G$ do
 not have the same length. By {\rm\cite{con}}, all subgroups of $G$ of order 2 are conjugate in $G$. Then $\delta(G)\geq |G:N_G(H)|=117>16$, a contradiction.

$(5)$ Let $G=S_z(2^q)$, where $q$ is an odd prime and $|G|=(2^{2q}+1)2^{2q}(2^q-1)$.

Arguing as in proof of Theorem~\ref{l1}, $G$ has a subgroup $H$ of order 2
all of whose maximal chains in $G$ do not have the same length and $|N_G(H)|\leq2^{2q}(2^q-1)$. Then $\delta(G)\geq|G:N_G(H)|\geq \frac{(2^{2q}+1)2^{2q}(2^q-1)}{2^{2q}(2^q-1)}
=2^{2q}+1\geq65>16$, also a contradiction.

All above arguments show that the counterexample of minimal order does not exist and so $G$ is solvable.\hfill\qed

\begin{lemma}\ \ \label{l22} Let $G$ be a non-solvable group. Then $\delta(G)=16$ if and only if $G\cong A_5$.
\end{lemma}

\demo We only need to prove the necessity part.

Since $G$ is non-solvable, there exists a subgroup $M$ of $G$ such that $M$ is a minimal non-solvable group. Then $M/\it\Phi$$(M)$ is a minimal
non-abelian simple group. Note that $\delta(M/\it\Phi$$(M))\leq\delta(M)\leq\delta(G)=16$, one must have $\delta(M/\it\Phi$$(M))=\delta(M)=\delta(G)=16$ by Lemma~\ref{l21}. Then arguing as
in proof of Lemma~\ref{l21}, one has $M/\it\Phi$$(M)\cong PSL_2(5)\cong PSL_2(4)\cong A_5$.

Claim $\it\Phi$$(M)=1$. Otherwise, assume $\it\Phi$$(M)\neq 1$. Since all maximal chains of the identity in $G$ do not have the same length, one has $\delta(G)\geq \delta(M/\it\Phi$$(M))+1$, a contradiction. Therefore, $\it\Phi$$(M)=1$.

It follows that $M\cong A_5$. If $M<G$. Since $\delta(M)=\delta(G)=16$, one has that $M$ is normal in $G$. Let $L$ be a subgroup of $G$ such that $M$ is maximal in $L$. Then $L/M\cong Z_p$ for some prime $p$.
One gets $L\cong S_5$ or $A_5\times Z_p$. However, it is easy to see that $\delta(S_5)>16$ and $\delta(A_5\times Z_p)>16$, a contradiction. Therefore, $G=M\cong A_5$. \hfill\qed

{\bf Proof of Theorem~\ref{th2}.} Combine Lemmas~\ref{l21} and ~\ref{l22} together, we complete the proof of Theorem~\ref{th2}.\hfill\qed

\section{Remarks}\label{s6}

\begin{remark}\ \ \label{r61} {\rm A group $G$ might not be solvable if $G$ has at least one minimal subgroup $H$ all of whose maximal chains in $G$ do not have the same length.
For example, let $G=SL_2(5)$. It is obvious that all maximal chains of every minimal subgroup of $G$ of odd order in $G$ have the same length. However, observe that
$Z(G)={\rm Z_2}<{\rm Z_4}<Q_8<SL_2(3)<G$ and $Z(G)={\rm Z_2}<{\rm Z_6}<SL_2(3)<G$ are two maximal chains of $Z(G)$ in $G$ having different length and $Z(G)$
is the unique subgroup of $G$ of order 2.}
\end{remark}

\begin{remark}\ \ \label{r62} {\rm A group $G$ might not be solvable if all maximal chains of any non-trivial subgroups of $G$ except minimal subgroups of order $p$
for a fixed prime divisor $p$ of $|G|$ have the same length. For example, it is easy to see that all maximal chains of non-trivial subgroups of $A_5$ except minimal subgroups of order 2
in $A_5$ have the same length, but $A_5$ is non-solvable.}
\end{remark}

\begin{remark}\ \ \label{r63} {\rm If all maximal chains of every second maximal subgroup of a group $G$ in $G$ have the same length, then $G$ might not be solvable. For example, let $G=PSL_2(17)$.
By {\rm\cite{DIC}}, it is easy to see that all maximal chains of every second maximal subgroup of $PSL_2(17)$ in $PSL_2(17)$ have the same length, but $PSL_2(17)$ is non-solvable.}
\end{remark}

\begin{remark}\ \ \label{r64} {\rm Assume that $G$ is solvable. If all maximal chains of every second maximal subgroup of $G$ in $G$ have the same length, then $G$ might not be supersolvable.
For example, let $G=GL_2(3)$. Let $H$ be any subgroup of $G$, then $H\in\{1,\,{\rm Z_2},\,{\rm Z_4},\,{\rm Z_8},\,{\rm Z_2\times Z_2},\, D_8,\, Q_8,\,\\SD_{16},\, SL_2(3),\,GL_2(3)\}$, where $SD_{16}=\langle a,\,b\mid\,a^8=b^2=1,\,b^{-1}ab=a^3\rangle$. It is easy to see that all maximal chains of every second maximal subgroup of $GL_2(3)$ in $GL_2(3)$ have the same length, but $GL_2(3)$ is non-supersolvable.}
\end{remark}

\end{document}